\documentclass{sigcomm-alternate}
\usepackage{amsfonts,latexsym,amsmath,amstext,amssymb,verbatim,epsfig,psfrag}
\usepackage{colordvi,pstcol}
\usepackage{graphicx}
\usepackage{psfrag}

\newcommand\II{\mathcal{I}}

\def\R{\mbox{I\hspace{-.15em}R}}

\def\Pb{\bold{P}}
\def\Ab{\bold{A}}

\def\Ib{\bold{I}}
\def\Jb{\bold{J}}

\newcommand\blfootnote[1]{%
  \begingroup
  \renewcommand\thefootnote{}\footnote{#1}%
  \addtocounter{footnote}{-1}%
  \endgroup
}
\newtheorem{definition}{Definition}
\newtheorem{remark}{Remark}
\newtheorem{theorem}{Theorem}
\newtheorem{lemma}{Lemma}
\newtheorem{property}{Property}

\begin{document}
\title{Convergence of the D-iteration algorithm: convergence rate and asynchronous distributed scheme}





\numberofauthors{3}
\author{
   \alignauthor{Dohy Hong}\\
   \affaddr{Alcatel-Lucent Bell Labs}\\
   \affaddr{Route de Villejust}\\
   \affaddr{91620 Nozay, France}\\
   \email{\normalsize dohy.hong@alcatel-lucent.com}
   \alignauthor {Fabien Mathieu}\\
   \affaddr{INRIA}\\
   \affaddr{avenue d'Italie}\\
   \affaddr{75014 Paris, France}\\
   \email{\normalsize fabien.mathieu@inria.fr}
   \alignauthor{G\'erard Burnside}\\
   \affaddr{Alcatel-Lucent Bell Labs}\\
   \affaddr{Route de Villejust}\\
   \affaddr{91620 Nozay, France}\\
   \email{\normalsize gerard.burnside@alcatel-lucent.com}
}

\maketitle

\begin{abstract}
In this paper, we define the general framework to describe the diffusion operators associated to a positive matrix.
We define the equations associated to diffusion operators and present some general properties of their state vectors.
We show how this can be applied to prove and improve the convergence of a fixed point problem associated to the matrix iteration scheme,
including for distributed computation framework.
The approach can be understood as a decomposition of the matrix-vector product operation in elementary operations at the vector entry level.
\end{abstract}
\category{G.1.3}{Mathematics of Computing}{Numerical Analysis}[Numerical Linear Algebra]
\category{G.1.0}{Mathematics of Computing}{Numerical Analysis}[Parallel algorithms]
\terms{Algorithms, Performance}

\keywords{Linear algebra, numerical computation, iteration, fixed point, eigenvector, distributed computation}
\blfootnote{*The work presented in this paper has been partially carried out at LINCS (www.lincs.fr).}
\begin{psfrags}
\section{Introduction}\label{sec:intro}
Today, we are living at the heart of the information age and we are facing
more and more data for which it becomes difficult to process using, say, traditional 
data processing strategies. In such a context,
iterative methods to solve large sparse linear systems have been gaining interests
in many different research areas and a large number of solutions/approaches have been studied.
For instance, iterative techniques gained a significant efficiency by exploiting
the sparsity of the information structure (decomposition, partition strategies etc).
However, if a near to optimal specific and direct solution can be built for a given
problem, it is obviously hard to have one solution that remains optimal for a large
classe of problems.
One of the most noticeable {\em generic} improvement was brought by iterative
methods combining the Krylov subspace and preconditioning approaches
\cite{Saad}, \cite{cg}, \cite{arnoldi}.

In this paper, we propose a new iterative algorithm based on a simple and intuitive fundamental
understanding of the linear equations as diffusions: 
we believe that this approach may bring significant improvement in a large classe of linear problems.
More precisely, we study the fixed point convergence problem in linear algebra exploiting the idea
of fluid diffusion associated to the D-iteration \cite{d-algo}.
This approach has been initially proposed to solve the PageRank equation \cite{serge}, \cite{dohy}.
The D-iteration solves $X\in\R^N$ of the equation:
$$
X = \Pb X + B
$$
where $\Pb$ is a non-negative matrix.
This includes in particular the case where $\Pb$ is of spectral radius unity (and $B=0$).

Solving a linear problem is a very well known theoretical problem and there are
a large number of methods that have been proposed, studied and explored.
For the general description of existing and/or alternative iteration methods,
one may refer to \cite{Golub1996}, \cite{Saad}, \cite{varga2009}, \cite{greenbaum}, \cite{stewart}.
In particular, 
for the power iteration method (solving $\Pb X = X$), the theory is very well known, for instance when $\Pb$ is 
associated to an irreducible transition matrix of a Markov chain, $X$ would
be its unique stationary distribution (cf. \cite{gantmacher2000theory}, \cite{berman1994nonnegative}).
The convergence of classical iterative schemes, such as Jacobi or Gauss-Seidel or successive over-relaxation method,
is also a very well known problem. One usual convergence condition is that (condition expressed for the equation
$\Ab X = B$, $\Ab$ may be chosen equal to $\Ib-\Pb$) $\Ab$ is strictly diagonal dominant.
With other approaches such as Krylov, conjugate gradient and variant methods exploiting the idea of projection
and residual minimization, better convergence can be obtained under more restrictive conditions on $\Ab$
(such as symmetric and positive definite).
If the power iteration (or the usual matrix-vector product) can be associated to
linear operations on the rows of the matrix (assuming an iteration is defined by the product matrix-vector), 
the diffusion approach consists in exploiting the columns of $\Pb$. 
The diffusion approach requires in particular that we have to define two state vectors $F$ (fluid)
and $H$ (history).
In this paper, we revisit the convergence condition of the D-iteration and we show that
its convergence can be obtained under a very wide condition on $\Pb$, with an upper bound on its convergence
rate when the spectral radius of $\Pb$ is strictly less than unity.
One of the main advantage of the diffusion approach is that its distance to the limit is explicitly known
and that the convergence speed gain can be intuitively and simply analyzed.\\

Concerning iterative methods on distributed computation framework, the usual approach is to
define a synchronization phase: for instance, in GPU computation \cite{gpu} or MapReduce framework \cite{mapreduce}, the
synchronization phase is an explicit step of the architecture. For those approaches, the convergence
condition is generally the same as for the sequential computation.
For an asynchronous distributed computation framework, the convergence may be harder to prove
or may simply fail:
\cite{mitra} is one of the first paper proving the convergence of the
asynchronous distributed computation of power iterations under pretty wide conditions.

The diffusion approach is well suited to an asynchronous distributed computation 
(cf. \cite{dist-test, part}), but its convergence had not been proved formally. 
In this paper, we formally prove the convergence of the diffusion approach (D-iteration) both for the
sequential computation (single processor) and for the distributed computation cases.

We compare the performance of the distributed computations based on the theoretical simulation
framework introduced in \cite{mitra} and show that for sparse large matrices, the case for
which the diffusion approach was initially designed, the computation
gain of the proposed method compared to the row-based asynchronous parallel computation
is significantly large, close to the theoretical optimal efficiency for large sparse matrices.\\

In Section \ref{sec:frame}, we define the notations and some general properties.
Section \ref{sec:diter} gives the formal framework of the D-iteration method
and some general results, in particular, the theoretical proofs of the
convergence for the sequential approach (single processor).
In Section \ref{sec:speed}, we show that for the most natural three variants of
the D-iteration, a simple upper bound of the convergence speed can be obtained
when $\Pb$ can be reduced in a simple multiplicative form.
Section \ref{sec:dist} gives the proof of the convergence of the distributed scheme.
Finally, Section \ref{sec:expe} reports and discusses some experimental results.
\section{General framework}\label{sec:frame}
We will use the following notations:
\begin{itemize}
\item $\Pb \in (\R^+)^{N\times N}$ a non-negative matrix;
\item $\Ib \in (\R^+)^{N\times N}$ the identity matrix;
\item $\Jb_i$ the matrix with all entries equal to zero except for
  the $i$-th diagonal term: $(\Jb_i)_{ii} = 1$;
\item $\Omega = \{1,..,N\}$;
\item $\II = \{i_1,i_2,..,i_n,...\}$ a sequence of nodes: $i_k \in \Omega$;
\item $\sigma_v : \R^N \to \R$ defined by $\sigma_v(X) = \sum_{i=1}^N v_i x_i$;
   if $V$ has no zero entry, we define the norm $|X|_v = \sum_{i=1}^N |v_i x_i|$;
\item $L_1$-norm: if $X\in \R^N$, $|X| = \sum_{i=1}^N |x_i|$;
\item $e$ the normalized unit column vector: $1/N \times (1,..,1)^T$;
\item $G = \{G_0, G_1, G_2,.., G_n, ..\}$ is a sequence of real vectors (in $\R^{N}$) such that $\sum_{n=0}^{\infty} |G_n| < \infty$;
\item $out_i$ is the number of non-zero entries of the $i$-th column of $\Pb$ (counts outgoing links).
\end{itemize}

We assume in this paper that $\Pb$ is non-negative for the sake of the simplicity.
We could generalize some results below when $\rho(\tilde{\Pb})\le 1$, where $\tilde{\Pb}$ is the matrix
where each component $p_{ij}$ is replaced by $|p_{ij}|$.

\subsection{Monotonicity}
We say that $\Pb$ is $\sigma_v$-decreasing if:
$$
\forall X\in(\R^+)^N, \sigma_v(\Pb X) \le \sigma_v(X).
$$

We define $\Pb^{\alpha} = (1-\alpha)\Ib + \alpha \Pb$.

Then, we have the following results:

\begin{property}\label{th:sigma}
$\sigma_v$-decreasing property is stable by composition of operators (matrix product).

If $\Pb$ is $\sigma_v$-decreasing, for all $\alpha \ge 0$, $\Pb^{\alpha}$ is $\sigma_v$-decreasing.

If $\Pb$ is $\sigma_v$-decreasing, for all $(\alpha,\alpha') \in (\R^+)^2$ such that $\alpha\le\alpha'$,
$\sigma_v(\Pb^{\alpha'}X) \le \sigma_v(\Pb^{\alpha}X)$.
\end{property}
\proof
The first point is obvious. The other points are based on the linearity of $\sigma_v$.

\subsection{Diffusion operators}
We define the $N$ diffusion operators associated to $\Pb$ by:
$$
\Pb_i = \Ib - \Jb_{i} + \Pb.\Jb_{i}
$$

\begin{property}\label{th:Pi}
If $\Pb$ is $\sigma_v$-decreasing, then the diffusion operators
$\Pb_i$ are $\sigma_v$-decreasing. Therefore, for $\alpha \ge 0$, 
$$\Pb^{\alpha}_i = (\Pb_i)^{\alpha} = \Ib + \alpha(\Pb-\Ib) \Jb_i$$ 
is $\sigma_v$-decreasing.
\end{property}
\proof
$\sigma_v(\Pb_i X) = \sigma_v(X) + \sigma_v(\Pb \Jb_i X) - \sigma_v(\Jb_i X)$ and we have
$\sigma_v(\Pb \Jb_i X) \le \sigma_v(\Jb_i X)$, therefore 
$\sigma_v(\Pb_i X) \le \sigma_v(X)$. The last point is the application of
Property \ref{th:sigma} to $\Pb_i$.

\section{Application to D-iteration}\label{sec:diter}
The D-iteration is defined by the couple $(\Pb,B) \in \R^{N\times N}\times \R^N$ and exploits two state vectors:
$H_n$ (history) and $F_n$ (residual fluid) based on the following iterative equations:
\begin{eqnarray}
F_0 &=& B\nonumber\\
F_n &=& \Pb_{i_n} F_{n-1}\label{eq:defF}
\end{eqnarray}
and
\begin{eqnarray}
H_0 &=& 0 =(0,..0)^T\nonumber\\
H_n &=& H_{n-1} + \Jb_{i_n} F_{n-1}.\label{eq:defH}
\end{eqnarray}

The D-iteration consists in updating the joint iterations on $(F_n, H_n)$.
Then, variant strategies may be applied depending on the choice of the sequence $\II$.
To recall the dependence on $\Pb,B$ and $\II$, we set:
$H_n(\Pb, B, \II) = H_n$. When the limit is well defined we will set
$H(\Pb, B, \II) = \lim_{n\to\infty}H_n(\Pb, B, \II)$.

We will consider two cases: if $\rho(\Pb)<1$ ($\rho$ is the spectral radius of $\Pb$), then 
we will see that $H_n(\Pb, B, \II)$ has a limit (denoted also $H_{\infty}$)
which is the unique solution of the equation (cf. Theorem \ref{theo:cv-m1}):
$$
X = \Pb X + B.
$$

If $\rho(\Pb) = 1$, then we will only consider in this paper the D-iteration with $B= \Pb.e - e$
(or any other vector for which $\sigma_v(B)=0$, with $V$ a positive left eigenvector of $\Pb$).

\subsection{General results}

\begin{theorem}[Fundamental diffusion equation]
We have:
$$H_{n+1}+F_{n+1} = H_n + F_n + \Pb (H_{n+1}-H_n)$$ and
\begin{eqnarray}
H_{n}+F_{n} &= \Pb H_n + F_0.\label{eq:fund}
\end{eqnarray}
\end{theorem}
\proof
The first equation is straightforward from the equations \eqref{eq:defF} and \eqref{eq:defH}.
The second one can be obtained by induction.\\

The first equation means that what we have (sum of $F$ and $H$) is what we had
before plus what's diffused by the increment of $H$.
The second equation means that
what you have is the initial value plus what you received from diffusion.

\begin{theorem}
We have:
\begin{eqnarray*}
H_n &=& \left(\Ib - \Jb_{i_n}(\Ib - \Pb)\right) H_{n-1} + \Jb_{i_n} B.
\end{eqnarray*}
\end{theorem}
\proof
We can rewrite the equation as
$$
H_n - H_{n-1} = \Jb_{i_n} \left( B - (\Ib - \Pb) H_{n-1}\right).
$$
Using \eqref{eq:defH}, we only need to check that $F_{n-1}$ is equal to
$B - (\Ib - \Pb) H_{n-1}$, which is exactly the equation \eqref{eq:fund}.

\subsection{Adding fluids}
Consider the D-iteration $H_n(\Pb,B,\II)$ on which we add $G$ (we will denote this by $H_n(\Pb,B,\II,G)$): before each diffusion
we add to $F_n$ the vector $G_n$.
This means that we modify $F_n$ and $H_n$ equations as follow:
$$
F_n = \Pb_{i_n} (F_{n-1} + G_{n-1})
$$
and
$$
H_n = H_{n-1} + \Jb_{i_n} (F_{n-1} + G_{n-1}).
$$

Then we have the following result:

\begin{theorem}
We have:
$$H_{n+1}+F_{n+1} = H_n + F_n + G_n + \Pb(H_{n+1}-H_n)$$ and
\begin{eqnarray}
H_{n}+F_{n} = \Pb H_n + F_0 + \sum_{i=0}^{n-1} G_i.\label{eq:fundg}
\end{eqnarray}
\end{theorem}

\subsection{Convergence}

We first show the convergence of $H_n(\Pb, B, \II)$ when $\rho(\Pb) < 1$.
This is a quite intuitive result and we only require that $\II$ is a
{\em fair} sequence.

\begin{definition}
A sequence $\II$ is fair if the number of occurrences of each $i\in\Omega$ is unbounded.
\end{definition}

\begin{remark}
If we have $i\in\Omega$ such that $(F_n)_i$ is equal to zero after finite steps $n_0$,
we don't need the fairness condition for the position $i$ (for the convergence).
\end{remark}

\begin{lemma}\label{lem:monof}
If $\Pb$ is irreducible and $\rho(\Pb)\le 1$, then there exists a (strictly) positive vector $V$ such that
$|F_n|_v$ is non-increasing. As a consequence, $F_n$ is convergent for
any given sequence $\II$.
\end{lemma}
\proof
Set $V$ the left positive eigenvector of $\Pb$ for $\rho(\Pb)$ 
($V$ is the left Perron vector, cf. \cite{gantmacher2000theory}, \cite{berman1994nonnegative}).\\
Set $j=i_{n+1}$ and $f=(F_n)_{i_{n+1}}$, then:
\begin{align*}
|F_{n+1}|_v =& \sum_{i\neq j} |(F_{n+1})_i v_i| + |(F_{n+1})_j v_j|\\
=& \sum_{i\neq j} |(F_n)_i + fp_{ij}| v_i + |fp_{jj} v_j|.
\end{align*}
Let's call $\Delta$ the set of nodes $i$ such that $(F_n)_i$ has a sign opposed to $f$.
Then,
\begin{align*}
|F_{n+1}|_v =& \sum_{i\neq j} (|(F_n)_i| + |f p_{ij}|) v_i + |f p_{jj}| v_j\\
    &+ \sum_{i\in \Delta} (|(F_n)_i + fp_{ij}| - |(F_n)_i| - |f p_{ij}|) v_i\\
=& |F_{n}|_v + \sum_{i\in \Delta} (|(F_n)_i + fp_{ij}| - |(F_n)_i| - |f p_{ij}|) v_i\\
\le & |F_{n}|_v.
\end{align*}
For the last inequality, we used $|x+y| \le |x|+|y|$.
Therefore, we have $|F_{n+1}|_v\le |F_n|_v$.
Therefore, $F_n$ is convergent. 

\begin{remark}
The above Lemma holds for any matrix $\Pb$, if there exists a positive vector $V$
such that for all $j$, $\sum_i (|\Pb_{ij}|\times v_i) \le v_j$.
\end{remark}

\begin{lemma}\label{lem:linH}
If $B= B_1 + B_2$, then for all $n$, we have
$$
H_n(\Pb,B,\II) = H_n(\Pb,B_1,\II) + H_n(\Pb,B_2,\II).
$$
\end{lemma}
\proof
The proof is straightforward using the linearity of $H_n$ w.r.t. $F_n$.

\subsubsection{Case $\rho(\Pb)<1$}

\begin{theorem}\label{theo:cv-m1}
If $\rho(\Pb) < 1$, for any fair sequence $\II$, $H_n(\Pb, B, \II)$ is convergent to the unique vector $X$
such that $X = \Pb X + B$.
\end{theorem}
\proof
Let's first assume that $B$ is non-negative, so that we only manipulate non-negative
quantities. By construction, $H_n$ is non-decreasing per entry.
From the equation \eqref{eq:fund}, we have:
$H_n = (\Ib - \Pb)^{-1} (B - F_n)$.
Hence, $H_n \le (\Ib - \Pb)^{-1} B$, therefore $H_n$ is convergent and because $\II$
is a fair sequence, necessarily, $F_n$ tends to zero.
Then, its limit satisfies the claimed equation.
Now, if $B$ has negative and positive terms, we can decompose $B$ as
$B^+ - B^-$ and apply the same argument for each component.

\begin{lemma}
If $\rho(\Pb)< 1$, then
for all $(\alpha, \beta) \in \R^2$,
$$\alpha H(\Pb,  B, \II) + \beta H(\Pb, B', \II) = H(\Pb, \alpha B + \beta B', \II).$$
\end{lemma}
\proof
We have:
\begin{align*}
H(\Pb, \alpha B + \beta B', \II) &= (\Ib-\Pb)^{-1}(\alpha B + \beta B')\\
&= \alpha (\Ib-\Pb)^{-1}B + \beta (\Ib-\Pb)^{-1}B'\\
&= \alpha H(\Pb,  B, \II) + \beta H(\Pb, B', \II).
\end{align*}

\begin{theorem}[Superposition]\label{theo:super}
Let $S = \sum_{n=0}^{\infty} G_n$ with $\sum_{n=0}^{\infty} |G_n| < \infty$.
If $\rho(\Pb)< 1$, then
$H(\Pb,B,\II,G) = H(\Pb,B+S,\II)$.
\end{theorem}
\proof
We have from the equation \eqref{eq:fundg}:
$H_n(\Pb,B,\II,G) = (\Ib-\Pb)^{-1}(B+\sum_{i=0}^{n-1} G_i - F_n)$.
Then, one can easily check that
the difference $|\sum_{i=n}^{\infty} G_i - F_n|$ tends to zero and
the equality holds.

\begin{theorem}[Monotonicity]\label{theo:mono-sup}
Let $S = \sum_{n=0}^{\infty} G_n$ with $\sum_{n=0}^{\infty} |G_n| < \infty$.
Let $S' = \sum_{n=0}^{\infty} G_n'$ with $\sum_{n=0}^{\infty} |G_n'| < \infty$
such that $G_n \le G_n'$.
If $\rho(\Pb)< 1$, then for all $n$,
$H_n(\Pb,B,\II,G) \le H_n(P,B,\II,G')$.
\end{theorem}
\proof
The proof is straightforward by induction using the iterative equations of $F_n$
and $H_n$ ($F_n\le F_n'$ and $H_n\le H_n'$).

\subsubsection{Case $\rho(\Pb)=1$}
In this case, we assume that the initial vector $B$ satisfies
$\sigma_v(B) = 0$, with $V$ the left eigenvector of $\Pb$ for
eigenvalue 1.

Let's first consider the case when $\Pb$ is irreducible.
We denote by $\Omega^+_n$ (resp. $\Omega^-_n$) the subset of $\Omega$ such that
$(F_n)_i \ge 0$ (resp. $(F_n)_i \le 0$).

\begin{lemma}\label{lem:iso}
The diffusion operators preserve $\sigma_v$, which means that:
$$
\sigma_v(\Pb_i X) = \sigma_v(X).
$$
\end{lemma}
\proof
$\sigma_v(\Pb_i X) = \sigma_v(X) + \sigma_v( (\Pb-\Ib) \Jb_i X)$.
And we have $\sigma_v( (\Pb-\Ib) \Jb_i X) = (V^T \Pb - V^T) \Jb_i X = 0$.\\

\begin{lemma}
For all $n$, $\sigma_v(F_n) = 0$.
\end{lemma}
\proof
This is a direct consequence of Lemma \ref{lem:iso} and $\sigma_v(F_0) = 0$ by
assumption.\\

\begin{lemma}
If at each diffusion step, $i_n\in \Omega^+_n$ (resp. $i_n\in \Omega^-_n$), then $\Omega^+_n$ 
(resp. $\Omega^-_n$) converges to $\Omega^+$ (resp. $\Omega^-$).
\end{lemma}
\proof
It is clear that the subset $\Omega^+_n$ is non-decreasing (we only add positive
quantities to each node). It is bounded by $\Omega$, hence convergent.\\

We denote by $\II^+$ (resp. $\II^-$) a fair sequence on $\Omega^+ = \cup_{n=1}^{\infty} \Omega^+_n$ 
(resp. $\Omega^-=\cup_{n=1}^{\infty} \Omega^-_n$),
then, we have the following result.

\begin{theorem}
If $\Pb$ is irreducible and $\II= \II^+$, then
$H_n(\Pb,\Pb.e-e,\II)$ is convergent to $X-e$ where $X$ is the real right eigenvector of $\Pb$
for the eigenvalue 1 with $\min_i x_i = 1/N$.
\end{theorem}
\proof
If there exists $T<\infty$ such that $\cup_{n=1}^{T} \Omega^+_n = \Omega$, 
then $\sigma_v(F^+_T) = \sigma_v(F^-_T) = 0$ and $H_n$ converged in finite time.
Otherwise, $\Omega^+$ is strictly included in $\Omega$. 
Let $\Pb^+$ be the restriction of $\Pb$ on $\Omega^+$: then $\rho(\Pb^+)<1$ (cf. \cite{berman1994nonnegative}).

At the limit, $H$ satisfies $(H+e) = \Pb (H+e)$ with $H+e$ a positive vector.
There is at least one coordinate $i$ on which the diffusion operator has been
never applied (with positive fluid), therefore $\min_i (H)_i = 0$.

\begin{remark}
If $\Pb$ is not irreducible, the diffusion on $\II^+$ may not converge. The counter example
is easy to be found.
\end{remark}


In order to prove the convergence for the sequence $\II^-$ without assuming
the irreducibility of $\Pb$, we consider a bit more general diffusion iterations as follows:
\begin{eqnarray}
F^{\alpha}_n &=& \Pb^{\alpha_n}_{i_n} F^{\alpha}_{n-1}\label{eq:defFg}
\end{eqnarray}
and
\begin{eqnarray}\label{eq:defHg}
H^{\alpha}_n &=& H^{\alpha}_{n-1} + \alpha_n \Jb_{i_n} F^{\alpha}_{n-1}
\end{eqnarray}
where $\alpha_n \ge 0$. If for all $n$, $\alpha_n=1$, we have the usual
diffusion iteration.

\begin{theorem}\label{th:HFgen}
$(F^{\alpha}_n, H^{\alpha}_n)$ satisfies:
\begin{align}\label{eq:defHpart}
H^{\alpha}_n + F^{\alpha}_n &= \Pb  H^{\alpha}_n + B.
\end{align}
\end{theorem}
\proof
The proof is the same as for the case $\alpha=1$, by induction
and using equations \eqref{eq:defFg} and \eqref{eq:defHg}.

\begin{theorem}
Assume we choose $F_0\ge 0$ and $H_0=0$. Then, $F^{\alpha}_n$ and
$H^{\alpha}_n$ are positive and $(H^{\alpha}_n)_i$ is an increasing function for all $i$. 

If $P$ is $\sigma_v$-decreasing, then
$\sigma_v(F^{\alpha}_n)$ is a decreasing function.
\end{theorem}
\proof
The proof is straightforward.

\begin{theorem}[Partial diffusion]\label{th:main}
If we build the two diffusion iterations $(F^{\alpha}_n, H^{\alpha}_n)$ and $(F^{\alpha'}_n,H^{\alpha'}_n)$
from the same initial vector $F_0$ ($H_0=0$) and for the same diffusion sequence $\II$,
if for all $n$, $0\le \alpha_n \le \alpha_n'\le 1$, then we have:
\begin{itemize}
\item $\sigma_v(F^{\alpha'}_n) \le \sigma_v(F^{\alpha}_n)$;
\item $H^{\alpha'}_n \ge H^{\alpha}_n$ (for each vector entry);
\item $H^{\alpha'}_n + F^{\alpha'}_n \ge H^{\alpha}_n + F^{\alpha}_n$.
\end{itemize}
\end{theorem}

\proof
The first inequality is a direct consequence of Property \ref{th:Pi} and Property \ref{th:sigma}.
For the second and third inequalities, we prove by induction: 
we have obviously $H^{\alpha'}_0 \ge H^{\alpha}_0$ and $H^{\alpha'}_0 + F^{\alpha'}_0 \ge H^{\alpha}_0 + F^{\alpha}_0$.
Assume we have, $H^{\alpha'}_n \ge H^{\alpha}_n$ and $H^{\alpha'}_{n} - H^{\alpha}_{n} \ge F^{\alpha}_{n}-F^{\alpha'}_n$. 
Then, from \eqref{eq:defHg}:
\begin{align*}
H^{\alpha'}_{n+1} &= H^{\alpha'}_{n} + \alpha_{n+1}' \Jb_{i_{n+1}} F^{\alpha'}_{n}\\
&\ge H^{\alpha'}_{n} + \alpha_{n+1} \Jb_{i_{n+1}} F^{\alpha'}_{n}
\end{align*}
and
\begin{align*}
H^{\alpha}_{n+1} &= H^{\alpha}_{n} + \alpha_{n+1} \Jb_{i_{n+1}} F^{\alpha}_{n}
\end{align*}

We first prove that:
\begin{align*}
H^{\alpha'}_{n} - H^{\alpha}_{n} &\ge \alpha_{n+1} \Jb_{i_{n+1}} (F^{\alpha}_{n}-F^{\alpha'}_n).
\end{align*}
For $i\neq i_{n+1}$, $(H^{\alpha'}_{n} - H^{\alpha}_{n})_i \ge 0$ and $(\Jb_{i_{n+1}} (F^{\alpha}_{n}-F^{\alpha'}_n))_{i} = 0$.
For $i = i_{n+1}$, we only need to handle the case $(F^{\alpha}_{n}-F^{\alpha'}_n)_{i_{n+1}} \ge 0$ and
we use the relation: $H^{\alpha'}_{n} - H^{\alpha}_{n} \ge F^{\alpha}_{n}-F^{\alpha'}_n$.
Hence, he have the inequality $H^{\alpha'}_{n+1} \ge H^{\alpha}_{n+1}$. 
Then, $H^{\alpha'}_{n+1} - H^{\alpha}_{n+1} \ge F^{\alpha}_{n+1}-F^{\alpha'}_{n+1}$
is straightforward using the equation \eqref{eq:defHpart} of Theorem \ref{th:HFgen}.

\begin{remark}
The power iteration $X_n = \Pb X_{n-1}$ can be described in the above
scheme  $(F^{\alpha}_n, H^{\alpha}_n)$
with $X_n = F^{\alpha}_{Nn}$
taking $F_0 = X_0$ and
if we apply the cyclic sequence $1,..,N$  ($i_n = n \mbox{ mod } N$) where $\alpha_{kN+i}\le 1$ is chosen
such that we diffuse exactly $(\Pb^k X_0)_i$
(such a value exists and is less than 1 because after the diffusion of nodes $1,..,i$
the residual fluid on $(i+1)$-th node can only  be increased).
\end{remark}
\begin{remark}
Note that the above theorem is valid without any condition on $\rho(P)$.
\end{remark}


\begin{theorem}\label{theo:neg-dif}
Assume we have a strictly positive left-eigenvector $V>0$ of $\Pb$.
If we choose any fair sequence of nodes $\II = \II^-$ (we only diffuse negative
fluids), then the diffusion applied on $(\Pb, \Pb e-e)$ converges to a unique $H_{\infty}$
such that $H_{\infty}+e$ is the right-eigenvector for eigenvalue
1 for $\Pb$, such that the maximum value within each strongly connected component of spectral radius one
is equal to $1/N$.
\end{theorem}
\proof
The diffusion from $\Pb.e-e$ can be decomposed as the difference of two
diffusion process $(F^{\alpha'}_n,H^{\alpha'}_n)$ and $(F^{\alpha}_n,H^{\alpha}_n)$ as follows:
we start with $F_0=e$. For the $N$ first diffusions, we choose $i_n=n$ and
\begin{itemize}
\item for $\Pb^{\alpha_n}_{i_n}$, $\alpha_n = 0$;
\item for $\Pb^{\alpha_n'}_{i_n}$, $\alpha_n'$ such that we diffuse exactly
$1/N$ from all nodes (such a value exists and is less than 1 because after the diffusion of nodes $1,..,i$
the residual fluid on $(i+1)$-th node can only be increased).
\end{itemize}
Then we have:
$F^{\alpha}_N = e$, $H^{\alpha}_N = 0$ and $F^{\alpha'}_N = \Pb.e$, $H^{\alpha'}_N = e$.
Then from the $(N+1)$-th diffusion, we apply exactly the same sequence with $\alpha_n = \alpha_n' = 1$.
For $n \ge N$, from Theorem \ref{th:main} (and linearity of $H$, Lemma \ref{lem:linH}, w.r.t. $F_n = F^{\alpha'}_n - F^{\alpha}_n$), 
we have $H_n + e = H^{\alpha'}_n - H^{\alpha}_n \ge 0$ and we have
$\sigma_v(F_n) = \sigma_v(F^{\alpha'}_n - F^{\alpha}_n) = 0$.
If we only diffuse negative fluids, this means that $H_n$ is a decreasing function (per entry).
Since we have $0 \le H_n + e \le e$, $H_n$ is convergent.
é
We proved above (Lemma \ref{lem:monof}) that $|F_n|_v = \sum_i |v_i\times (F_n)_i|$ is a decreasing function.
The convergence of $H_n$ implies of course the convergence to zero of $F_n$ on the coordinate
$v_i > 0$ and we have $H+e = \Pb (H+e)$.
Now we still have to prove that $H+e \neq 0$. In fact, there exists at least one
entry such that $(H)_i = 0$: take $i\in \Omega^+_n$ (here we require that $\Omega^+_n$ is only for
entries strictly positive). If the iteration stops in finite
time, the node in the last non empty $\Omega^+_n$ satisfies $(H)_i = 0$. Otherwise, there exists
at least one node such that it stays always strictly positive, for which we have $(H)_i = 0$.
Hence, $\max_i (H+e)_i = 1/N$. The argument here applies for each strongly
connected component (SCC) of spectral radius one (SCC(1)). 
For the uniqueness, we can decompose $\Pb$ in three parts (cf. \cite{berman1994nonnegative}): 
\begin{itemize}
\item those who point to SCC(1)s and are not pointed to by any SCC(1)s: for them
  we have  $(H)_i=0$;
\item those in SCC(1)s: for each SCC(1), $H$ restricted on a SCC(1) is equal to the unique eigenvector
  of SCC(1) (Perron eigenvector) such that the maximum is equal to $1/N$;
\item those who are pointed to by SCC(1)s and are not pointing to any SCC(1)s: this case corresponds
  to the previous section (diffusion on $\rho < 1$) for which the initial condition results
  from $\Pb e-e$ and the fluids received from SCC(1)s, on which we have uniqueness.
\end{itemize}

In order to justify clearly the above decomposition, we argue as follows:
because of the existence of a strictly positive left-eigenvector of $\Pb$ (for 1),
we can decompose $\Omega$ as in Figure \ref{fig:decomp}: 

\begin{figure}[htbp]
\centering
\includegraphics[width=8cm]{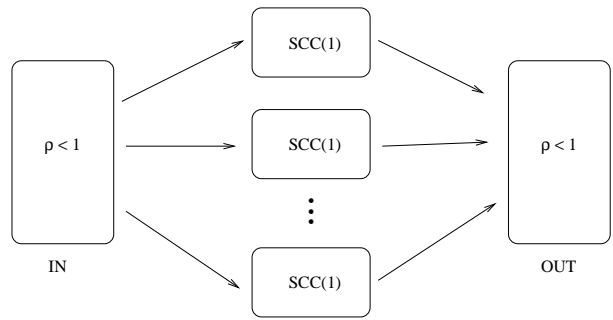}
\caption{Decomposition.}
\label{fig:decomp}
\end{figure}

where the restriction of $\Pb$ on the subsets IN and OUT has a spectral radius
strictly less then 1.

The diffusion of nodes in IN is independent of the rest and its diffusion limit
is $(H+e)_{i\in IN} = 0$ as shown in the following Lemma:

\begin{lemma}
If $\rho(\Pb) < 1$, and if we apply only negative diffusion from the initial condition
$F_0=P.e-e$, then the limit is $H = -e$.
\end{lemma}
\proof
Since $\rho(\Pb)<1$, the convergence of $H$ is clear. We need to justify why diffusing
only negative fluid, we do not end up with some strictly positive coordinates.
From Theorem \ref{th:main}, we have $H_n+e \ge 0$. Then from the equation \eqref{eq:fund},
one can easily get: $\sigma_v(F_n) = (\rho(\Pb)-1) \sigma_v( H_n + e)$, which means that
$F_n$ has always some negative components. At the limit, we can not have any negative
fluids and $\sigma_v(F_{\infty}) = (\rho(\Pb)-1) \sigma_v( H_{\infty} + e)$, which implies
$H_{\infty}=-e$ and $F_{\infty} = 0$.

From the decomposition of Figure \ref{fig:decomp}, once we have that $(H)_i + 1/N = 0$ on the nodes in IN,
the diffusion of each SCC(1) is exactly the independent diffusion of $\Pb$ restricted to each
SCC(1) component, which converges to the unique eigenvector for 1.
Finally, the nodes of OUT converges to a unique solution, since again we have $\rho(OUT)<1$ and
we can apply Theorem \ref{theo:super}. This ends the proof of the theorem.

\begin{remark}
If $V$ has some coordinates equal to zero, the above convergence works only on
$i$ such that $v_i>0$ and we can have $(F_{\infty})_i > 0$ on $i$ such that $v_i=0$.
Let $\Omega^{+\infty}$ be the set coordinates such that $(F_{\infty})_i > 0$.
Then the spectral radius of $P$ restricted to $\Omega^{+\infty}$ may be one, in which
case, there is no convergence. 
\end{remark}


\begin{remark}
If we mix the diffusion of positive and negative fluids, there is no guarantee
that the D-iteration algorithm converges. For instance, with a snake-configuration counter example,
we may oscillate (take for instance, $\Pb$ associated to a five nodes graph: 
$1\to 2$ and $1\to 3$; $2\to 4$; $3\to 5$; $4\to 1$; $5\to 1$).
\end{remark}


\begin{theorem}[Monotonicity 2]\label{theo:mono-sup2}
Let $S = \sum_{n=0}^{\infty} G_n$ with $\sum_{n=0}^{\infty} |G_n| < \infty$.
Let $S' = \sum_{n=0}^{\infty} G_n'$ with $\sum_{n=0}^{\infty} |G_n'| < \infty$.
If for all $n$, $\sum_{i=0}^{n} G_i \le \sum_{i=0}^{n} G_i'$,
then for all $n$, $H_n(\Pb,B,\II,G) \le H_n(\Pb,B,\II,G')$.
\end{theorem}
\proof
This is a particular case of the partial diffusion result (Theorem \ref{th:main}):
$H_n(\Pb,B,\II,G)$ can be seen as $H_n(\Pb,B+G',\II)$ where at step $n$, 
$\sum_{i=n+1}^{\infty} G_i' + \sum_{i=0}^n (G_i'-G_i)$ has been blocked.
$H_n(\Pb,B,\II,G')$ can be seen as $H_n(\Pb,B+G',\II)$ where at step $n$, 
$\sum_{i=n+1}^{\infty} G_i'$ has been blocked.

\section{Convergence speed}\label{sec:speed}

For the sake of the simplicity, we only considered here the case of
sub-stochastic matrices of the form: $\Pb= d.\Pb_s$ where
$0<d<1$ and $\Pb_s$ is a transition matrix.

We consider the following choice of sequence:
\begin{itemize}
\item CYC: $i_n = n \mod N$;
\item MAX: $i_n =\arg\max_i (F_n)_i$;
\item COST: $i_n =\arg\max_i \left((F_n)_i/out_i\right)$.
\end{itemize}

Let's first consider the theoretical cost of the iteration as the number of times
diffusions are applied.

\begin{theorem}
The convergence speed of CYC and MAX is at least $d^{\lfloor l/N\rfloor}$, where
$l$ is the number time diffusions are applied.
\end{theorem}
\proof
CYC: from $F_n$, after $N$ diffusions, we diffused at least $|F_n|$, therefore,
$|F_{n+N}| \le d.|F_n|$. \\
MAX: if we order $(F_n)_i$: $(F_n)_1 \ge (F_n)_2 \ge ... (F_n)_N$, we start by diffusing $(F_n)_1$, and
at the $k$-th diffusion we take a fluid at least equal to the $k$-th term $(F_n)_k$. As for CYC, we diffused at least $|F_n|$
after $N$ diffusions.\\

Now let's consider the theoretical cost of the iteration as the number of times
a non-zero element of $P$ is used.

\begin{theorem}
The convergence speed of COST is at least $d^{\lfloor l/L\rfloor}$,
where $L= \sum_i out_i$ and $l$ is the number of times
a non-zero element of $\Pb$ is used.
\end{theorem}
\proof
We just need to count as above the amount of fluid we move while using $L$ links of $P$.
At step $n$, assume we order $(F_n)_i/out_i: (F_n)_1/out_1 \ge (F_n)_2/out_2 \ge .. \ge (F_n)_N/out_N$.
At $n+1$, we diffuse exactly $(F_n)_1$ which costs $out_1$.
Assume we have diffused $D_{k-1}$ with $L_{k-1}$ operations. Then,
At step $n+k$, we diffuse $d_k$ with cost $l_k$ such that $d_k/l_k \ge (F_n)_k/out_k$.
Hence, we have $D_k/L_k \ge (\sum_{i=1}^k (F_n)_i)/(\sum_{i=1}^k out_i) \ge |F_n|/L$.
If we stop counting when $L_k$ just exceeds $L$ (say $k_0$), then,
$D_{k_0} \ge |F_n|/L\times L_{k_0} \ge |F_n|$.\\

\begin{remark}
An interesting heuristic bound is $d/(2-d)$ assuming both the graph and the sequence are random.
\end{remark}

\section{Convergence of distributed algorithm}\label{sec:dist}

\subsection{Case $\rho < 1$}
Here we first consider the case $\rho(\Pb)<1$.
We detail the proof of the convergence for $K=2$ when the computation is done on two
virtual processors, that we call PID1 and PID2. We assume a fixed partition of
$\Omega = \Omega_1 \cup \Omega_2$ and the corresponding decomposition of $\Pb$ as:

$$
P = 
\begin{pmatrix}
\Pb_{11} & \Pb_{12}\\
\Pb_{21} & \Pb_{22}
\end{pmatrix}
$$

The diffusion state variables of two PIDs are denoted by:
$(F^1,H^1)$, $(F^2,H^2)$. For a single PID sequential system, we note
$H = ([H]^1,[H]^2)$.

To prove the convergence of the distributed architecture to the proper eigenvector,
we consider the following decomposition: we assume that $PID_1$ computes the D-iteration $H(\Pb_{11}, B_1, \II^1)$
with additional inputs from $PID_2$ at times $T_1,..T_n,...$.
During $\tau_n = T_n - T_{n-1}$ (we set $T_0=0$), $PID_1$ applies the diffusion sequence
$\II^1_n$.
We assume without any loss of generality that $PID_2$ computes the D-iteration based
on the sequence $\II^2$ ($\II^2_n$ during $T_n' - T_{n-1}'$) and sends at $T_n'$ ($n\ge 1$) the quantity 
$$
\Pb_{12}(H^2_{T_n'}-H^2_{T_{n-1}'}) {\buildrel \text{def}\over =} \Delta H^2_n.
$$
Let's first assume that this quantity defines the additional input at $T_n$ by $PID_1$ (one way).
The information transmission delay is then defined by $T_n - T_n'$.
We will see that we don't need that the information from $PID_2$ arrives in the order for the
proper convergence. The only requirement is that $T_n - T_n'$ and $T_n'-T_{n-1}'$ are bounded.

\begin{figure}[htbp]
\psfrag{PID1}{$PID_1$}
\psfrag{PID2}{$PID_2$}
\psfrag{I1}{$\II^1_1$}
\psfrag{I2}{$\II^1_2$}
\psfrag{Inp}{$\II^1_{n+1}$}
\psfrag{J1}{$\II^2_1$}
\psfrag{J2}{$\II^2_2$}
\psfrag{Jnp}{$\II^2_{n+1}$}
\psfrag{T1}{$T_1$}
\psfrag{T2}{$T_2$}
\psfrag{Tn}{$T_n$}
\psfrag{Tnprime}{$T_n'$}
\psfrag{Tnp}{$T_{n+1}$}
\psfrag{deltan}{$\Pb_{12}\Delta H^2_n$}
\centering
\includegraphics[width=8cm]{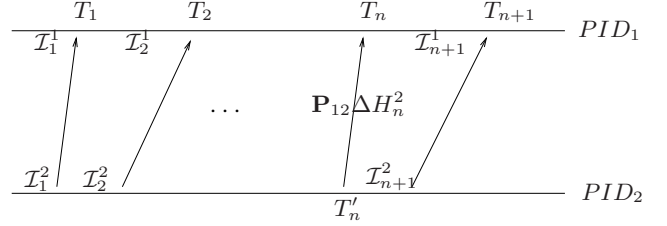}
\caption{Data exchange: from PID2 to PID1.}
\label{fig:ex}
\end{figure}

\noindent
If $n$ is such that $T_k \le n < T_{k+1}$, we have (upper bound):
\begin{align*}
H^1_n =& H_{T_1}(P_{11}, B_1, \II^1_1)\\
       & + H_{\tau_2}(P_{11}, F^1_{T_1}+\Delta H^2_1, \II^1_2)\\
       & ...\\
       & + H_{\tau_k}(P_{11}, F^1_{T_{k-1}}+\Delta H^2_{k-1},\II^1_k)\\
       & + H_{n-T_k}(P_{11}, F^1_{T_k}+\Delta H^2_k,\II^1_{k+1})\\
      =& (\Ib-\Pb_{11})^{-1}(B_1 + \Pb_{12}H^2_{T_k} - F^1_n).
\end{align*}
Hence,
\begin{align}\label{eq:ine}
 H^1_n     \le& (\Ib-\Pb_{11})^{-1}(B_1 + \Pb_{12}H^2_{n}).
\end{align}
In the above, we used the equality
$$
H_{\tau_i}(\Pb_{11}, F^1_{T_{i-1}}+\Delta H^2_{i-1},\II)
= (\Ib-\Pb_{11})^{-1} (F^1_{T_{i-1}}+\Delta H^2_{i-1} - F^1_{T_{i}}).
$$
The above result shows that only the cumulated fluids matter for PID1.\\

Let's now use the notation of the $G$ (fluid addition): if $n=T_i$, we set $G^2_n = \Delta H^2_i$,
otherwise $G^2_n = 0$. We assume $S=\sum_{i=1}^{\infty} G^2_i = \Pb_{12} H^2_{\infty} < \infty$.
Then $H^1_n = H_n(\Pb_{11}, B_1, \II^1, G^2)$ (we define $G^1_n$ in the same way).
From Theorem \ref{theo:super}, we have
$H^1_{\infty} = H(\Pb_{11}, B_1+G, \II^1)$.\\




Now consider the global system. 
Let's denote by $(T_n',T_n)$ the sending time and the receiving time of the
$n$-th fluid transmission exchange between PIDs (in both directions). Because the different PIDs may
have different computation speed, we may have that at $T_n$ the number of iterations
done by the PID that receives the fluid is less than the number of iterations
done by the sender at time $T_n'$. For the sake of the simplicity, we still denote
by $F_{T_n}$ (same for $H$ and for $T_n'$) the value of the $F$ at time $T_n$.

\begin{lemma}[Impact of the delay]\label{lem:delay}
Consider two distributed computation systems $S(1), S(2)$ that have exactly the same $T_n'$
and that differs only by $T_n$ such that for all $n$, $T_n(1) \le T_n(2)$.
Then, we have: $H^1_n(1) \ge H^1_n(2)$ and $H^2_n(1) \ge H^2_n(2)$.
\end{lemma}
\proof
This is a direct consequence of Theorem \ref{theo:mono-sup2}.

\begin{lemma}\label{lem:borne}
At each step, he have $H^1_n \le [H]^1_n$ and $H^2_n \le [H]^2_n$.
Therefore, the distributed computation scheme converges.
\end{lemma}
\proof
From Lemma \ref{lem:delay}, we first have the monotonicity when
reducing the delay up to zero delay $T_n = T_n'$. Then, again thanks to Theorem \ref{theo:mono-sup2},
one can show that adding more fluid exchanges, we can only increase the history vector.
The system with zero delay and exchange at each step corresponds to the sequential system $H$.
Therefore, we have  $H^1_n \le [H]^1_n\le [H]^1_{\infty}$
and $H^2_n \le [H]^2_n\le [H]^2_{\infty}$, and $H^1_n$ and $H^2_n$ are convergent.

\begin{theorem}
If each PID applies a fair sequence for diffusions and if
for all $n$, $|T_n - T_n'| \le T$ and $|T_n'-T_{n-1}'| \le T$,
then the distributed computation converges to the sequential
single PID computation $H$.
\end{theorem}

\proof
From the equation \eqref{eq:ine}, we still have:
$H^1_{n} \le (\Ib-\Pb_{11})^{-1}(B_1 + \Pb_{12}H^2_{\infty})$ and
$H^2_{n} \le (\Ib-\Pb_{22})^{-1}(B_2 + \Pb_{21}H^1_{\infty})$.

If for all $n$, $|T_n - T_n'| < T$ and $|T_n'-T_{n-1}'| < T$
we have (lower bound):
$$
H^1_{n+2T} \ge (\Ib-\Pb_{11})^{-1}(B_1 + \Pb_{12}H^2_{n} - F^1_{n+2T})
$$
and
$$
H^2_{n+2T} \ge (\Ib-\Pb_{22})^{-1}(B_2 + \Pb_{21}H^1_{n} - F^2_{n+2T})
$$

We already know that both converges (non-decreasing and upper bounded
by $[H]^1_{\infty}$ and $[H]^2_{\infty}$), therefore $F^1_{n+2T}$
and $F^2_{n+2T}$ goes to zero and at the limit we have:

$$
H^1_{\infty} = (\Ib-\Pb_{11})^{-1}(B_1 + \Pb_{12}H^2_{\infty})
$$
and
$$
H^2_{\infty} = (\Ib-\Pb_{22})^{-1}(B_2 + \Pb_{21}H^1_{\infty}).
$$

Therefore, we have $H^1_{\infty} = [H]^1_{\infty}$ and
$H^2_{\infty} = [H]^2_{\infty}$ (because $H$ is the unique solution
satisfying the above two equations).\\

Now, the generalization of this proof to any $K>1$ follows exactly
the same ideas.

\subsection{Case $\rho = 1$}
Here, we assume that there exists $V$ a strictly positive
left eigenvector of $\Pb$ and that the initial vector $B$ satisfies
$\sigma_v(B) = 0$.
Let's consider the limit of the diffusion on $\Omega^-$.

We use the same notation than in the previous section.
Since we only apply the diffusion on nodes having negative fluids,
the history vector $H$ is for each entry a negative decreasing function.

In the previous case with $\rho<1$, we could assume that $\II^1$ and $\II^2$ are
predefined and independent of the fluid exchanges.
The first difficulty here is that a priori we can no more do such an assumption,
since the sign of $(F)_i$ depends on the past fluid exchanges.

To overcome this difficulty and apply the previous results, 
we do the following {\em trick} which consists in
the decomposition of the diffusion in two steps: selection of a node and diffusion test.
We assume given two fair sequences $\II^1$ and $\II^2$ on $\Omega_1$ and $\Omega_2$
(fair means here that every coordinate will be candidate for diffusion an infinite
number of times) and when a node is selected, we only diffuse if its sign is
negative.

Then, we have the following result:
\begin{theorem}[Monotonicity extended]\label{theo:mono-sup3}
Let $G_n$ and $G_n'$ two sequences of non-positive vectors, such that
for all $n$, $\sum_{i=1}^{n} G_i \le \sum_{i=1}^{n} G_i'$. Let $\II$ be a fair sequence.
If we only diffuse the negative fluids,
then for all $n$, $H_n(\Pb,B,\II,G) \le H_n(\Pb,B,\II,G')$.
\end{theorem}
\proof
The argument is exactly the one used in Theorem \ref{theo:mono-sup2}.\\

Then, the results of Lemma \ref{lem:delay} and Lemma \ref{lem:borne} still hold with inversed inequality.

\begin{theorem}\label{theo:gen-dist}
If each PID applies a fair sequence for diffusions and if
for all $n$, $|T_n - T_n'| \le T$ and $|T_n'-T_{n-1}'| \le T$,
then the distributed computation on $\Omega^-$ converges to the sequential
single PID computation $H$ obtained by the diffusion of negative fluid.
\end{theorem}
\proof
We have the lower bounds:
$[H]^1_{\infty} \le [H]^1_n \le H^1_n$ and $[H]^2_{\infty} \le [H]^2_n \le H^2_n$,
therefore $H^1_n$ and $H^2_n$ are convergent.


We know that the limit of the distributed computation satisfies:
\begin{align}\label{eq:ine-1}
(\Ib-\Pb_{11}) H^1_{\infty}     =& B_1 + \Pb_{12}H^2_{\infty}\\
(\Ib-\Pb_{22}) H^2_{\infty}     =& B_2 + \Pb_{21}H^1_{\infty}.
\end{align}

Let $Z = ((H^1_{\infty})^T, (H^2_{\infty})^T)^T - H$.
If $\rho(\Pb_{11}) < 1$ and $\rho(\Pb_{22}) < 1$, then we have the uniqueness of
the solution and $Z=0$.
If $\rho(\Pb_{11}) = 1$ and $\rho(\Pb_{22}) = 1$, then since $\rho(\Pb)=1$ and
since we assumed the existence of strictly positive left-eigenvector for 1, 
we can not have any positive terms linking $\Omega_1$ and $\Omega_2$
(cf. \cite{berman1994nonnegative}), 
hence we have two independent systems for which the uniqueness was proved 
(Theorem \ref{theo:neg-dif}).
Finally, if $\rho(\Pb_{11}) = 1$ and $\rho(\Pb_{22}) < 1$, then we must have
$\Pb_{12}=0$, and for the same reason, we still have $Z=0$.

\section{Experimentation}\label{sec:expe}
\subsection{Uniform graph with $N=128$}
We first report the experiment results by reproducing the hypothetical parallel
computer simulations as defined in \cite{mitra} (Table 1). We reproduced a random (symmetrical)
graph of 128 nodes containing 1652 links, 23 of which are self-loops. 
As in \cite{mitra}, we assumed that a processor spends 0 cycles to read and write its local
memory, $T_r$ and $T_w$ cycles, respectively, to read from and write into the shared memory,
$T_m$ cycles for one multiplication and $T_a$ cycles for one addition
(below, we took $T_r=4$, $T_w=2$, $T_m=T_a=1$).
The degrees
of the nodes vary from 14 to 38, with a mean of 25.6 and standard deviation of
4.9. In Table \ref{tab:1}, we show the results of:
\begin{itemize}
\item sPI-R: synchronized power iteration per row (i.e. Jacobi iteration);
\item aPI-R: asynchronized power iteration per row (method evaluated in \cite{mitra});
\item sPI-C: synchronized power iteration per column;
\item sPI-Cr: synchronized power iteration per column assuming identical weight $p_{ij}$
  for each column $j$ (as in PageRank equation);
\item DI$+$COST: asynchronized diffusion applying COST from the initial vector $\Pb.e-e$
 (cost of the computation of $\Pb.e-e$ is included).
\end{itemize}

We see that the results of sPI-R and aPI-R are very close to those observed in
Table 1 of \cite{mitra}. With sPI-C, we see a gain for $K=2,4,8$: that's because
the fluid exchange to other PIDs aggregate the results on the diffusions of
$N/K$ nodes to $N/K*(K-1)$ nodes. However, there is a higher penalty when $K$ becomes
close to $N$. This can be explained considering the limit case of $K=N$: 
whereas one coordinate updates on a row having $n_r$ non-zero values
requires about $2\times n_r$ additions and multiplications
followed by one writing cycle, the diffusion would require $2\times n_c$ additions and
multiplications (for a column having $n_c$ non-zero values) followed by $n_c$
updates (read, addition and write) of shared memory.
The result for sPI-Cr shows the gain when exploiting the homogeneity of the column weight
(one multiplication required per diffusion), which can be done only when column based
operations are used (as for the diffusion approach).

\begin{table}
\begin{center}
\begin{tabular}{|l|ccccc|}
\hline
K   & sPI-R & aPI-R & sPI-C & sPI-Cr& DI$+$COST\\
\hline
1   & 65620 & 45934 & 65620 & 34090 & 18183\\
2   & 66020 & 45472 & 37500 & 21470 & 20047\\
4   & 41620 & 28504 & 23180 & 15030 & 17347\\
8   & 24000 & 16590 & 16000 & 11770 & 13711\\
16  & 13060 & 9142  & 11940 & 9780  & 9952\\
32  & 7880  & 5056  & 9980  & 8720  & 6469\\
64  & 4120  & 3056  & 8210  & 7550  & 4264\\
128 & 2260  & 2256  & 7215  & 6845  & 3128\\
\hline
\end{tabular}\caption{Uniform graph $N=128$.}\label{tab:1}
\end{center}
\end{table}

\begin{figure}[htbp]
\centering
\includegraphics[angle=-90,width=8cm]{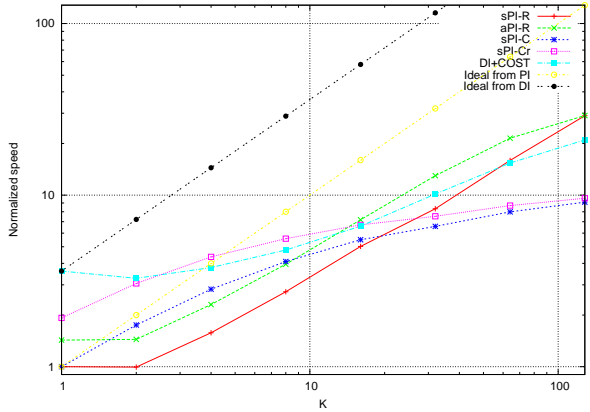}
\caption{Normalized speeds.}
\label{fig:mitra}
\end{figure}

Figure \ref{fig:mitra} shows the results of Table \ref{tab:1} in terms of
the convergence speeds normalized by the value of sPI-R for $K=1$.
Ideal-PI and Ideal-DI are the ideal curve ($y=x$) starting from $K=1$
of PI and DI.

\subsection{Web graph}
In this section, we used the 
web graph imported from the dataset \verb+uk-2007-05+ \verb+@1000000+
(available on \cite{webgraphit}) which has
41,247,159 links on $10^6$ nodes.
The aim is to compare the theoretical computation cost for a large sparse
matrix such as the one associated to a web graph, the case for which the
diffusion approach was initially designed.
The $N=1000$ case is obtained considering the first 1000 nodes of
the above graph.
Here the convergence is on the PageRank equation (eigenvector with
damping factor of $0.85$). The results are shown on Figure \ref{fig:mitra-1000} for
$N=1000$ and Figure \ref{fig:mitra-1000000} for $N=10^6$.

\begin{figure}[htbp]
\centering
\includegraphics[angle=-90,width=8cm]{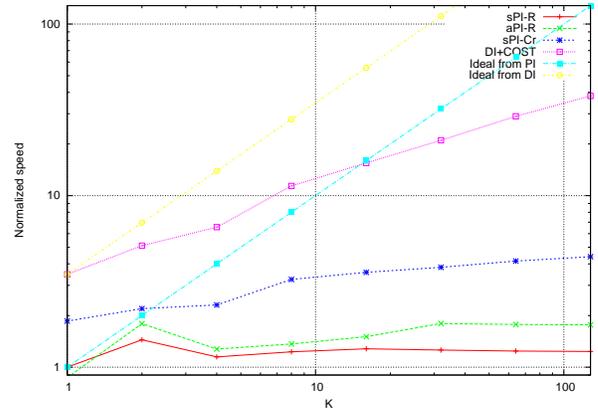}
\caption{Normalized speeds: $N=1000$.}
\label{fig:mitra-1000}
\end{figure}

In Figure \ref{fig:mitra-1000}, the benefit of the column based diffusion approach
is much more significant and even more significant for $N=10^6$ (Figure \ref{fig:mitra-1000000}). 
In Figure \ref{fig:mitra-1000000}, we added the performance of a dynamical partition approach
(cf. \cite{part} for more details: the dynamical partition used here consists
roughly in observing the convergence speed in logscale of each PID based on its remaining fluid
quantity and transferring 10\% of nodes that is managed by the slowest PID to the fastest PID
when the difference is higher than 50\%) in order
to check/validate the property claimed in \cite{part}: we can see that applying
the cost assumption/model of \ref{fig:mitra-1000}, a simple dynamical partition strategy on the
diffusion method leads to a performance that is close to the optimal efficiency (close to
the ideal curve).
Note that the cost of the
partition updates/modifications is partially included here 
(unity cost per node exchanged on both sides of PIDs).

\begin{figure}[htbp]
\centering
\includegraphics[angle=-90,width=8cm]{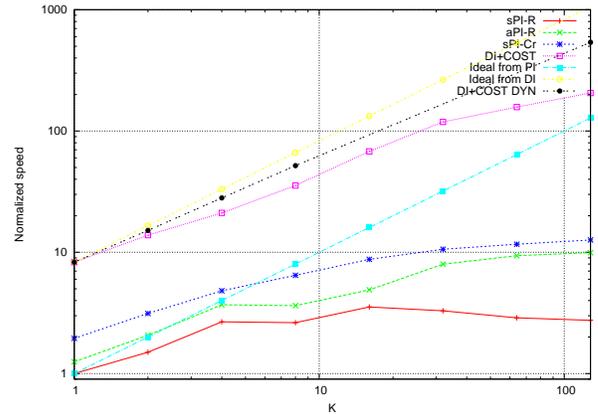}
\caption{Normalized speeds: $N=10^6$.}
\label{fig:mitra-1000000}
\end{figure}

\section*{Acknowledgments}
The authors are very grateful to Philippe Jacquet for his very valuable comments
and suggestions.

\section{Conclusion}\label{sec:conclusion}
In this paper, we addressed the formal proof of convergence of different
D-iteration schemes, including the case of the distributed computation
and the upper bounds of the convergence speeds.
We used the theoretical formal approach of \cite{mitra} to evaluated the
theoretical computation cost and showed the significant gain of our
diffusion approach for the computation of the eigenvector of large sparse
matrices.

\end{psfrags}
\bibliographystyle{abbrv}
\bibliography{sigproc}

\end{document}